\documentclass[10pt]{amsart}
\usepackage{amssymb,amsfonts,lineno,amsthm,amscd,stmaryrd}
\usepackage[leqno]{amsmath}
\usepackage[numbers]{natbib}

\theoremstyle{plain}
\newtheorem{thm}{Theorem}
\newtheorem{cor}[thm]{Corollary}
\newtheorem{lem}[thm]{Lemma}
\newtheorem{prop}[thm]{Proposition}
\theoremstyle{definition}
\newtheorem{definition}[thm]{Definition}
\newtheorem{remark}[thm]{Remark}
\newtheorem{example}[thm]{Example}

\numberwithin{thm}{section}
\numberwithin{equation}{section}

\newcommand{\EQ}[1]{\eqref{eq:#1}} 
\newcommand{\LEM}[1]{Lemma~\ref{lem:#1}}    
    
\newcommand{\THM}[1]{Theorem~\ref{thm:#1}}  
\newcommand{\REM}[1]{Remark~\ref{rem:#1}}  
\newcommand{\PROP}[1]{Proposition~\ref{prop:#1}}  
\newcommand{\COR}[1]{Corollary~\ref{cor:#1}} 
\newcommand{\SEC}[1]{Section~\ref{sec:#1}}

\newcounter{hypo}

\newcommand{\ep}{\varepsilon}

\newcommand{\Ellip}{\Lambda} 
\newcommand{\ellip}{\lambda}

\newcommand{\scon}{b}      

\newcommand{\limitset}{\mathcal{S}}

\DeclareMathOperator{\trace}{trace}

\newcommand{\Falpha}[1]{F_\alpha \left[ #1 \right]}

\newcommand{\puccisub}[2]{\mathcal{P}^-_{#1,#2}}
\newcommand{\Puccisub}[2]{\mathcal{P}^+_{#1,#2}}
\newcommand{\PucciSub}{\Puccisub{\ellip}{\Ellip}}
\newcommand{\pucciSub}{\puccisub{\ellip}{\Ellip}}
\newcommand{\pucci}{\mathcal{P}^-}
\newcommand{\Pucci}{\mathcal{P}^+}

\newcommand{\R}{\ensuremath{\mathbb{R}}}

\newcommand{\Sy}{\ensuremath{\mathbb{S}^n}}

\newcommand{\normx}[1]{\left\| #1 \right\|_*}   
\newcommand{\decayseta}{\mathcal{C}}   

\newcommand{\rescale}[1]{\mathcal{T}_{#1}}

\begin{document}
\title[Asymptotics for fully nonlinear parabolic equations]{Long-time asymptotics for fully nonlinear homogeneous parabolic equations}
\author{Scott N. Armstrong}
\address{Department of Mathematics,
University of California, Berkeley, CA 94720.}
\address{Department of Mathematics,  
Louisiana State University, Baton Rouge, LA, 70803}
\email{armstrong@math.lsu.edu}
\author{Maxim Trokhimtchouk}
\address{Department of Mathematics,
University of California, Berkeley, CA 94720.}
\email{trokhim@math.berkeley.edu}
\date{\today}

\begin{abstract}
We study the long-time asymptotics of solutions of the uniformly parabolic equation
\begin{equation*}
u_t + F(D^2u) = 0 \quad \mbox{in} \ \R^n\times \R_+,
\end{equation*}
for a positively homogeneous operator $F$, subject to the initial condition $u(x,0) = g(x)$, under the assumption that $g$ does not change sign and possesses sufficient decay at infinity. We prove the existence of a unique positive solution $\Phi^+$ and negative solution $\Phi^-$, which satisfy the self-similarity relations
\begin{equation*}
\Phi^\pm (x,t) = \lambda^{\alpha^\pm} \Phi^\pm ( \lambda^{1/2} x, \lambda t ).
\end{equation*}
We prove that the rescaled limit of the solution of the Cauchy problem with nonnegative (nonpositive) initial data converges to $\Phi^+$ ($\Phi^-$) locally uniformly in $\R^n \times \R_+$. The anomalous exponents $\alpha^+$ and $\alpha^-$ are identified as the principal half-eigenvalues of a certain elliptic operator associated to $F$ in $\R^n$. 
\end{abstract}

\maketitle


\section{Introduction and main results}

The connection between the scaling invariance of the mathematical expressions for certain physical laws and the asymptotic behavior of physical phenomena is of fundamental importance to the study of mechanics. In this work, we have in mind the study of self-similar solutions of  diffusion equations and their relationship to the long-time behavior of solutions to the Cauchy problem. A classical example is the heat equation
\begin{equation}\label{eq:heatequation}
u_t - \Delta u = 0 \quad \mbox{in} \ \R^n \times \R_+,
\end{equation}
which is invariant with respect to any scaling $(x,t) \mapsto \left(\sigma^{1/2} x , \sigma t\right)$, with $\sigma > 0$. It is well-known that a solution $u(x,t)$ of equation \EQ{heatequation} with nonnegative and integrable initial data will converge in a certain sense to a multiple of the Gaussian kernel
\begin{equation*}
\Phi(x,t) := (4\pi t)^{-\frac{n}{2}} e^{-\frac{|x|^2}{4t}}.
\end{equation*}
Precisely, the rescaled solutions $u^\sigma (x,t) := \sigma^{n/2} u\left( \sigma^{1/2} x, \sigma t\right)$ will converge locally uniformly to $C \Phi(x,t)$, where the constant $C$ is given by
\begin{equation*}
C = \int_{\R^n} u(x,0) \, dx.
\end{equation*} 
The Gaussian kernel $\Phi$ satisfies the relation
\begin{equation*}
\Phi(x,t) = \sigma^{\frac{n}{2}} \Phi\left( \sigma^{1/2} x, \sigma t\right), 
\end{equation*} 
and for this reason $\Phi$ is called a \emph{self-similar solution} of \EQ{heatequation}.

\medskip

Another important example of a diffusion equation is the generalized porous medium equation
\begin{equation}\label{eq:porous-medium-equation}
u_t - \Delta \Psi(u) = 0 \quad \mbox{in} \ \R^n \times \R_+,
\end{equation}
where $\Psi$ is an increasing function of $u$. The study of self-similar solutions of \EQ{porous-medium-equation} and the long-time asymptotics of solutions to the Cauchy problem is well-developed, particularly in the case that $\Psi(u) = u^m$ for $m>1$. It originated with the celebrated work of Barenblatt, Zel'dovich, and Kompaneets, and was later developed Friedman, Kamin, V\'azquez, and others (see V\'azquez \cite{Vazquez:Book} for a well-written introduction to the subject as well as a comprehensive list of references).

\medskip

While the literature on self-similar solutions and asymptotics of diffusion equations is vast, relatively little is known in the case of the fully nonlinear parabolic equation
\begin{equation} \label{eq:fully-nonlinear-parabolic}
u_t + F(D^2u) = 0 \quad \mbox{in} \ \R^n \times \R_+.
\end{equation}
Here $F$ is a positively homogeneous, uniformly elliptic operator. Important examples of \EQ{fully-nonlinear-parabolic} include the parabolic Bellman and parabolic Isaacs equations, which arise in the theory of stochastic optimal control and stochastic differential games, respectively. In the article, we show that equation \EQ{fully-nonlinear-parabolic} possesses a unique positive and negative self-similar solution, and that these self-similar solutions characterize the long-time asymptotic behavior of solutions of the Cauchy problem with initial data which does not change sign.

\begin{thm}\label{thm:existence-self-similar-solutions}
Assume that $F$ satisfies \EQ{Felliptic} and \EQ{Fhomogeneous}, below. Then there exist unique constants $\alpha^+=\alpha^+(F)>0$ and $\alpha^-=\alpha^-(F) > 0$, for which the uniformly parabolic equation \EQ{fully-nonlinear-parabolic} possesses solutions $\Phi^+ > 0$ and $\Phi^- < 0$, satisfying the relations
\begin{equation*}
\Phi^+(x,t) = \sigma^{\alpha^+} \Phi^+\left(\sigma^{1/2} x, \sigma t\right) \quad \mbox{and} \quad \Phi^-(x,t) = \sigma^{\alpha^-} \Phi^-\left(\sigma^{1/2} x, \sigma t\right)
\end{equation*}
for all $\sigma > 0$, and such that for some constants $C, a> 0$,
\begin{equation*}
\Phi^+(x,1) , - \Phi^-(x,1) \leq C \exp( - a|x|^2) \quad \mbox{for all} \ x \in \R^n.
\end{equation*}
Moreover, the solutions $\Phi^+$ and $\Phi^-$ are unique up to multiplication by a positive constant.
\end{thm}

We call the numbers $\alpha^+$ and $\alpha^-$ the \emph{positive} and \emph{negative anomalous exponents} of the operator $F$, respectively. In contrast to the case of a linear parabolic equation, $\alpha^\pm \neq n/2$ and $\alpha^+ \neq \alpha^-$, in general. We will see in \SEC{existence} that they are the principal half-eigenvalues of an elliptic operator in $\R^n$. The functions $\Phi^+(\cdot,1)$ and $\Phi^-(\cdot,1)$ are the corresponding principal half-eigenfunctions.

Our second main result characterizes the long-time behavior of solutions of equation \EQ{fully-nonlinear-parabolic} with initial data which does not change sign and exhibits sufficient decay at infinity.

\begin{thm} \label{thm:convergence-self-similar-solutions}
Let the hypotheses of \THM{existence-self-similar-solutions} be in force, and consider a viscosity solution $u \in C(\R^n \times \R_+)$ of the uniformly parabolic equation \EQ{fully-nonlinear-parabolic} such that $|u(x,0)| \leq C_0\exp (-B|x|^2)$ for some constants $B,C_0>0$. If $u(\cdot,0) \geq 0$ and $u(\cdot,0) \not\equiv 0$, then there exists a constant $C^* > 0$ such that the rescaled solutions given by
\begin{equation*}
u^\sigma (x,t):= \sigma^{\alpha^+} u \left( \sigma^{1/2} x, \sigma t\right), \quad (x,t) \in \R^n \times \R_+
\end{equation*}
converge locally uniformly in $\R^n \times \R_+$ as $\sigma \to \infty$ to the function $C^*\Phi^+$. Likewise, if $u(\cdot,0) \leq 0$, $u(\cdot,0) \not \equiv 0$, then there exists a constant $C^* > 0$ such that the functions 
\begin{equation*}
u^\sigma (x,t):= \sigma^{\alpha^-} u \left( \sigma^{1/2} x, \sigma t\right), \quad (x,t) \in \R^n \times \R_+
\end{equation*}
converge locally uniformly to $C^* \Phi^-$ as $\sigma \to \infty$.
\end{thm}

\medskip

The \emph{Barenblatt equation of elasto-plastic filtration}
\begin{equation}\label{eq:Barenblatt}
u_t - \max\left\{ \frac{\Delta u}{1-\gamma}, \frac{\Delta u}{1+\gamma} \right\} = 0 \quad \mbox{in} \ \R^n \times \R_+,
\end{equation}
where $0 < \gamma < 1$, is a particular example of a fully nonlinear parabolic equation of the form \EQ{fully-nonlinear-parabolic}. It arises in filtration theory as a model for an elastic fluid flowing through an irreversibly deformable elasto-plastic porous medium (see Barenblatt, Entov, and Ryzhik \cite{Barenblatt:Book:1990}). Kamin, Peletier and V\'azquez \cite{Kamin:1991} proved Theorems \ref{thm:existence-self-similar-solutions} and \ref{thm:convergence-self-similar-solutions} for equation \EQ{Barenblatt}, using different methods from those employed in this paper. Their argument makes use of the rotational invariance of the equation to reduce the problem to an ODE. The authors then employ a shooting argument in the phase plane to demonstrate the existence of self-similar solutions of equation \EQ{Barenblatt}. Their proof of the asymptotic convergence of rescaled solutions of the Cauchy problem to a multiple of the self-similar solution relies on a careful analysis of the self-similar solution at infinity. These estimates also rely on ODE theory in a crucial way, and in particular some clever applications of l'H\^opital's rule.

\medskip

Our proof of \THM{convergence-self-similar-solutions} employs the ``four-step method" of Kamin and V\'azquez \cite{Kamin:1988} in a new setting. The idea is to squeeze the rescaled functions $u^\sigma$ between two multiples of the self-similar solution $\Phi^\pm$, and to show that the gap between them must vanish asymptotically. Similar to \cite{Kamin:1991}, the primary difficulty in carrying out this analysis is controlling the behavior of the solutions as $t^{-1/2} |x| \to \infty$. This is overcome by the construction of a special comparison function.

\medskip

Self-similar solutions of equation \EQ{Barenblatt} have also been considered by Caffarelli and Stefanelli \cite{Caffarelli:2008}, and the dependence of the anomalous exponents on the parameter $\gamma$ has been examined by Goldenfeld, Martin, Oono, and Liu \cite{Goldenfeld:1990} as well as Aronson and V\'azquez \cite{Aronson:1995}. We also wish to mention the interesting work of Hulshof and V\'azquez \cite{Hulshof:1996}, who studied the long-time asymptotic behavior of viscosity solutions of an equation similar to \EQ{Barenblatt}, but with $\Delta(u^m)$ replacing both instances of $\Delta u$.

As we were completing a final revision of this paper for publication, we became aware of a new preprint by Meneses and Quaas \cite{Quaas:preprint} in which the existence assertions of \THM{existence-self-similar-solutions} are also obtained.

\medskip

In \SEC{notation} we state our hypotheses and recall the definition of a viscosity solution. \SEC{existence} contains the proof of \THM{existence-self-similar-solutions}, and in \SEC{asymptotics} we prove \THM{convergence-self-similar-solutions}.


\section{Notation and Hypotheses}\label{sec:notation}

Throughout this paper, we denote $(0,\infty)$ by $\R_+$. The set of $n$-by-$n$ real symmetric matrices is $\Sy$. For $M \in \Sy$ and $0 < \lambda \leq \Lambda$, define the operators
\begin{equation*}
\PucciSub(M) := \sup_{A\in \llbracket\lambda,\Lambda\rrbracket} \left[ - \trace(AM) \right] \quad \mbox{and} \quad \pucciSub (M) := \inf_{A\in \llbracket\lambda,\Lambda\rrbracket} \left[ - \trace(AM) \right],
\end{equation*}
where $\llbracket\lambda,\Lambda\rrbracket \subseteq \Sy$ is the set of positive definite matrices with eigenvalues contained in the interval $[ \lambda, \Lambda]$. The nonlinear operators $\PucciSub$ and $\pucciSub$ are called the \emph{Pucci maximal} and \emph{minimal operators}, respectively. For ease of notation, we will often drop the subscripts and write $\Pucci$ and $\pucci$. A convenient way to write the Pucci extremal operators is
\begin{equation}\label{eq:pucci-nice-form}
\Pucci(M) = -\lambda \sum_{\mu_j > 0} \mu_j - \Lambda \sum_{\mu_j < 0} \mu_j \quad \mbox{and} \quad \pucci(M) = -\Lambda \sum_{\mu_j > 0} \mu_j - \lambda \sum_{\mu_j < 0} \mu_j,
\end{equation}
where $\mu_1, \ldots, \mu_n$ are the eigenvalues of $M$.

\medskip

We require our operator $F:\Sy \rightarrow \R$ to be uniformly elliptic in the sense that there exist constants $0<\lambda \leq \Lambda$, such that
\begin{equation} \label{eq:Felliptic}
\puccisub{\lambda}{\Lambda}(M-N) \leq F(M) - F(N) \leq \Puccisub{\lambda}{\Lambda}(M-N) \quad \mbox{for all} \quad M,N\in \Sy.
\end{equation}
We also require $F$ to be positively homogeneous of order one:
\begin{equation}\label{eq:Fhomogeneous}
F(\eta M) = \eta F(M)\quad \mbox{for all} \quad M \in \Sy, \  \eta \geq 0.
\end{equation}

\begin{remark}\label{rem:its-all-the-same}
If $F$ satisfies hypotheses \EQ{Felliptic} and \EQ{Fhomogeneous}, then the operator $\tilde{F}(M):= - F( -M)$ satisfies \EQ{Felliptic} and \EQ{Fhomogeneous} as well. This observation will simplify the proofs of our main results. For example, the existence of the negative self-similar solution for $F$ can be deduced from the existence of the positive self-similar solution for $\tilde{F}$, and it is clear that $\alpha^-(F) = \alpha^+(\tilde{F})$. Likewise, to prove \THM{convergence-self-similar-solutions}, we need only to assume nonnegative initial data and show convergence to a multiple of the positive self-similar solution, as the last statement follows from the former applied to $\tilde{F}$.
\end{remark}

Every differential equation and differential inequality in this paper is assumed to be satisfied in the viscosity sense. An introduction to the theory of viscosity solutions can be found in Crandall, Ishii, and Lions \cite{UsersGuide}. For the convenience of the reader, we now recall the definition of a viscosity solution.

\begin{definition}
Let $G$ be an elliptic nonlinear operator, and $\Omega$ an open set in $\R^n \times \R_+$. A function $u\in C(\Omega)$ is a \emph{viscosity subsolution (supersolution)} of the parabolic equation
\begin{equation}\label{eq:defviscsol}
u_t + G(D^2u,Du,u,x) = 0 \quad \mbox{in} \ \Omega
\end{equation}
if whenever $x_0\in \Omega$ and $\varphi\in C^2(\Omega)$ are such that
\begin{equation*} 
x \mapsto u(x) - \varphi(x) \quad \mbox{has a local maximum (minimum) at} \quad x_0,
\end{equation*}
we have
\begin{equation*}
u_t + G(D^2\varphi(x_0),D\varphi(x_0),u(x_0),x_0) \leq \ (\geq)\ 0.
\end{equation*}
We say that $u$ is a \emph{viscosity solution} of \EQ{defviscsol} if it is both a viscosity subsolution and supersolution of \EQ{defviscsol}. Viscosity (sub/super)solutions of elliptic equations are defined analogously. The precise meaning of the differential inequality
\begin{equation*}
u_t + G(D^2u,Du,u,x) \leq \ (\geq) \ 0 \quad \mbox{in} \ \Omega
\end{equation*}
is that $u$ is a viscosity subsolution (supersolution) of \EQ{defviscsol}.
\end{definition}

Several times in this article, we will make use of the strong maximum principle and estimates for viscosity solutions of uniformly parabolic equations. For these facts we refer to Wang \cite{Wang:1992a} and Crandall, Kocan, and \'Swi\c{e}ch \cite{Crandall:2000}. Analogous results for uniformly elliptic equations can be found, for example, in Caffarelli and Cabre \cite{Caffarelli:Book}, Trudinger \cite{Trudinger:1988}, and Winter \cite{Winter:2009}.


\section{Existence of self-similar solutions} \label{sec:existence}

In this section we prove \THM{existence-self-similar-solutions}. To motivate our approach, suppose that there exists a solution $\Phi > 0$ of equation \EQ{fully-nonlinear-parabolic} and an exponent $\alpha > 0$ which satisfy the relation
\begin{equation}\label{eq:self-similarity-relation}
\Phi(x,t)  = \sigma^{\alpha} \Phi\left( \sigma^{1/2} x, \sigma t\right) \quad \mbox{for every} \ (x,t) \in Q, \ \sigma > 0.
\end{equation}
For $t=1$ this reads
\begin{equation} \label{eq:motivation-existence}
\Phi(x,1)  = \sigma^{\alpha} \Phi\left( \sigma^{1/2} x, \sigma \right).
\end{equation}
Formally differentiate \EQ{motivation-existence} with respect to $\sigma$ at $\sigma =1$ to discover that
\begin{equation} \label{eq:motivation-existence-1}
\Phi_t \left( x, 1\right) + \alpha \Phi(x,1) + \frac{1}{2} x \cdot D \Phi(x,1) = 0.
\end{equation}
Inserting \EQ{motivation-existence-1} into the PDE \EQ{fully-nonlinear-parabolic} and rearranging, we derive the equation
\begin{equation*}
F\left( D^2 \Phi(x,1) \right) - \frac{1}{2} x \cdot D\Phi(x,1) = \alpha \Phi(x,1) \quad \mbox{for all} \ x \in \R^n.
\end{equation*}
Defining $\varphi(y) := \Phi(y,1)$, we see that the pair $(\alpha,\varphi)$ is a solution of the elliptic eigenvalue problem
\begin{equation} \label{eq:elliptic-eigenvalue-Rn}
F\left( D^2\varphi \right) - \frac{1}{2} y \cdot D\varphi = \alpha \varphi \quad \mbox{in} \ \R^n.
\end{equation}

To prove \THM{existence-self-similar-solutions}, we will reverse the calculation above. Studying the elliptic eigenvalue problem \EQ{elliptic-eigenvalue-Rn}, we will show that there exists $\alpha > 0$ and $\varphi > 0$ which satisfy \EQ{elliptic-eigenvalue-Rn} and such that $\varphi$ decays at a suitable rate as $|y| \to \infty$. We will then define $\Phi$ by 
\begin{equation*}
\Phi(x,t) := t^{-\alpha} \varphi\left( t^{-1/2} x \right),
\end{equation*}
and check that $\Phi$ is a solution of \EQ{fully-nonlinear-parabolic} which also satisfies \EQ{self-similarity-relation}.

\medskip

Another way to derive \EQ{elliptic-eigenvalue-Rn} is to consider the \emph{continuous rescaling} of \EQ{fully-nonlinear-parabolic} given by the change of variables
\begin{equation*}
y = t^{-1/2} x, \ s = \log t, \ \varphi(y,s) = t^{\alpha} \Phi(x,t),
\end{equation*}
and then to look for a stationary solution $\varphi(y,s) = \varphi(y)$. Thus the appearance of \EQ{elliptic-eigenvalue-Rn} in our context is similar to the emergence of the so-called Fokker-Planck equation in the study of self-similar solutions of the porous medium equation (see V\'azquez \cite[Section 18.4]{Vazquez:Book}). An elliptic equation similar to \EQ{elliptic-eigenvalue-Rn} also arises in the study of the self-similar solutions of the semilinear parabolic equation
\begin{equation*}
u_t - \Delta u = | u|^{p-1} u \quad \mbox{in} \ \R^n \times \R_+.
\end{equation*}
See the work of Peletier, Terman, and Weissler \cite{Peletier:1986} and Haraux and Weissler \cite{Haraux:1982}.

\medskip

The theory of principal eigenvalues of fully nonlinear operators goes back to the work of Lions \cite{Lions:1983d}, who used stochastic methods to study the Bellman equation in a bounded domain. Recently, several authors, including Birindelli and Demengel \cite{Birindelli:2006,Birindelli:2007}, Ishii and Yoshimura \cite{Ishii:preprint}, and Quaas and Sirakov \cite{Quaas:2008} have studied the principal eigenvalues of more general fully nonlinear operators in bounded domains (see also \cite{Armstrong:2009}). Our methods in this section are similar to those in these works, although extra complications arise from the unboundedness of $\R^n$. The special form of \EQ{elliptic-eigenvalue-Rn}, particularly the gradient term, provides us with enough compactness to overcome these obstacles.

\begin{lem}\label{lem:SDMP}
Let $\alpha \in \R$, $r>0$, and assume that $u\in C(\R^n\backslash B_r)$ satifies
\begin{equation}\label{eq:SDMP}
\left\{ \begin{aligned}
& \pucciSub(D^2u) - \frac{1}{2} y\cdot Du \leq \alpha u &  \mbox{in} & \ \R^n\backslash B_r, \\
& u \leq 0 & \mbox{on} & \ \partial B_r,
\end{aligned} \right.
\end{equation}
and $u(y) \leq C e^{-a|y|^p}$ for some constants $C,a> 0$ and $p>1$. Then there exists a constant $R=R(\alpha,\Lambda)>0$ such that $r\geq R$ implies that $u\leq 0$.
\end{lem}
\proof
Set $R:= 2(\alpha + \Lambda+1)$. A routine calculation verifies that the function $\varphi(y):= \exp(-|y|)$ satisfies
\begin{equation}\label{eq:SDMP-1}
\pucci \left(D^2\varphi(y) \right) - \frac{1}{2} y \cdot D\varphi(y) \geq \left( -\Lambda + |y|/2 \right) \varphi(y) \geq (\alpha + 1)\varphi  \quad \mbox{in} \ \{ |y| \geq R \}.
\end{equation}
Suppose on the contrary that there exists $|y_0| > r$ such that $u(y_0) > 0$. Let 
\begin{equation*}
\eta := \inf \left\{ \rho > 0 : \rho \varphi \geq u \ \mbox{in} \ \R^n \backslash B_r \right\}.
\end{equation*}
Then $\eta > 0$, $\eta \varphi \geq u$ in $\R^n \backslash B_r$, and owing to the faster decay rate of $u$, there exists $|y_1| > r$ such that $\eta \varphi(y_1) = u(y_1)$. In particular, the map
\begin{equation*}
y \mapsto u(y) - \eta\varphi(y) \quad \mbox{has a local maximum at} \ y = y_1.
\end{equation*}
Recalling that $u$ is a viscosity solution of \EQ{SDMP}, we see that
\begin{equation*}
\pucci\left(\eta D^2\varphi(y_1) \right) - \eta \frac{1}{2} y_1\cdot D\varphi(y_1) \leq \alpha u(y_1) = \eta \alpha \varphi(y_1) < \eta (\alpha + 1) \varphi(y_1),
\end{equation*}
a contradiction to \EQ{SDMP-1}.
\qed

\medskip

The next lemma is an analogue of \cite[Theorem 3.3]{Armstrong:2009}, adapted to our setting. The result goes back to observations due by Berestycki, Nirenberg and Varadhan \cite{Berestycki:1994}, and earlier versions have been developed by Quaas and Sirakov \cite{Quaas:2008}, and Ishii and Yoshimura \cite{Ishii:preprint} to study the principal half-eigenvalues of fully nonlinear elliptic operators in bounded domains. It is an important comparison tool which is essential to our analysis in this section.

\begin{lem}\label{lem:parabolic-HCP}
Assume that $u,v,f \in C(\R^n)$, $f\geq 0$, and $\alpha \in \R$ satisfy
\begin{equation*}
F(D^2u) - \frac{1}{2} y\cdot Du - \alpha u \leq f \leq F(D^2v) - \frac{1}{2} y\cdot Dv - \alpha v \quad \mbox{in} \ \R^n.
\end{equation*}
Suppose also that $v>0$ and $u(y) \leq C\exp(-a |y|^p)$ for some $p>1$ and every $y\in \R^n$. Then either $u \leq v$ in $\R^n$, or $u \equiv t v$ for some constant $t > 1$.
\end{lem}
\proof
Let $R = R(\alpha,\Lambda)$ be as in \LEM{SDMP}. For $s\geq 1$, define $w_s:= u - sv$. Then
\begin{align*}
\puccisub{\lambda}{\Lambda}(D^2w_s) - \frac{1}{2} y \cdot Dw_s - \alpha w_s \leq f - sf \leq 0 \quad \mbox{in} \ \R^n.
\end{align*}
For any $s \geq 1$ so large that $w_s < 0$ on $\bar{B}_R$, \LEM{SDMP} implies that $w_s \leq 0$ on $\R^n \backslash B_R$. Hopf's lemma then implies that $w_s < 0$ in $\R^n$. Define
\begin{equation*}
t := \inf\left\{ s \geq 1: w_s < 0 \ \mbox{in} \ \R^n \right\}.
\end{equation*}
Then $t\geq 1$ and $w_t \leq 0$. If $t = 1$, then $u \leq v$ in $\R^n$, and we have nothing left to show. Suppose that $t > 1$. We claim that $w_t \equiv 0$. If not, then by Hopf's lemma, $w_t < 0$. Select $\delta > 0$ so small that $t - \delta > 1$ and $w_{t-\delta} < 0$ on $\bar{B}_R$. Then $w_{t-\delta} < 0$ on $\R^n$, as we argued above. This contradicts the definition of $t$, completing the proof.
\qed

\begin{cor}\label{cor:parabolic-HCP}
Suppose that $u,v \in C(\R^n)$ and $\alpha \in \R$ satisfy
\begin{equation*}
F(D^2u) - \frac{1}{2} y\cdot Du - \alpha u \leq 0 \leq F(D^2v) - \frac{1}{2} y\cdot Dv - \alpha v \quad \mbox{in} \ \R^n.
\end{equation*}
Suppose also that $v>0$ and $u(y) \leq C\exp(-a |y|^p)$ for some $p>1$ and every $y\in \R^n$. Then either $u \leq 0$ in $\R^n$, or $u \equiv t v$ for some constant $t > 0$.
\end{cor}
\proof
Suppose that $u(y_0) > 0$ for some $y_0\in \R^n$. Using the homogeneity of $F$ and by multiplying $u$ by a large positive constant, if necessary, we may assume that $u(y_0) > v(y_0)$. According to \LEM{parabolic-HCP}, there exists $t > 0$ such that $v \equiv t u$. 
\qed

\medskip

In our analysis, an important role with be played by the functions $\left\{ \phi^a \right\}_{a>0}$, which we define by
\begin{equation*}
\phi^a(y) := \exp\left(-a|y|^2\right).
\end{equation*}

\begin{lem}\label{lem:rugged-calculations}
For $a=1/4\Lambda$, the function $\phi^a$ satisfies
\begin{equation}
\left( \frac{n\lambda}{2\Lambda}\right) \phi^a \leq \pucciSub(D^2\phi^a) - \frac{1}{2} y\cdot D\phi^a \leq \left( \frac{(n-1)\lambda}{2\Lambda} + \frac{1}{2} \right) \phi^a \quad \mbox{in} \ \R^n.
\end{equation}
Likewise, for $a=1/4\lambda$, the function $\phi^a$ satisfies
\begin{equation}
\left( \frac{(n-1)\Lambda}{2\lambda} + \frac{1}{2} \right) \phi^a \leq \PucciSub(D^2\phi^a) - \frac{1}{2} y\cdot D \phi^a \leq \left( \frac{n\Lambda}{2\lambda}\right) \phi^a  \quad \mbox{in} \ \R^n.
\end{equation}
\end{lem}
\proof
For any $a> 0$, the Hessian of $\phi^a$ is given by
\begin{equation*}
D^2\phi^a(y) = \phi^a(y) \left( 4a^2y \otimes y - 2aI \right),
\end{equation*}
and has eigenvalues $(4a^2|y|^2 - 2a)\phi^a(y)$ with multiplicity $1$ and $-2a\phi^a(y)$ with multiplicity $n-1$. Hence
\begin{equation*}
\pucciSub(D^2\phi^a(y)) = \begin{cases}
a \phi^a(y) (2 \lambda n - 4a\lambda |y|^2), & \mbox{for all} \  |y| \leq (2a)^{-\frac{1}{2}} , \\
a \phi^a(y) ( 2 \lambda (n-1) + 2\Lambda - 4a\Lambda |y|^2), & \mbox{for all} \  |y| > (2a)^{-\frac{1}{2}}.
\end{cases}
\end{equation*}
Therefore, for each $a> 0$,
\begin{equation}\label{eq:pucci-bound-alpha-1}
\pucciSub(D^2\phi^a) - \frac{1}{2}y\cdot D\phi^a = a\left(2\lambda n - 4a\lambda |y|^2 + |y|^2\right) \phi^a\quad \mbox{for} \  |y| \leq (2a)^{-\frac{1}{2}},
\end{equation}
and 
\begin{equation}\label{eq:pucci-bound-alpha-2}
\pucciSub(D^2\phi^a) - \frac{1}{2}y\cdot D\phi^a = a\left(2\lambda (n-1) +2\Lambda - 4a\Lambda |y|^2 + |y|^2\right) \phi^a \quad \mbox{for} \  |y| \geq (2a)^{-\frac{1}{2}}.
\end{equation}
In particular, for every $a\leq (4\Lambda)^{-1}$, we have
\begin{equation}\label{eq:pucci-bound-alpha}
\pucciSub(D^2\phi^a) - \frac{1}{2}y\cdot D\phi^a \geq (2a\lambda n) \phi^a \quad \mbox{in} \ \R^n.
\end{equation}
For $a=(4\Lambda)^{-1}$ we carefully check that
\begin{equation}\label{eq:pucci-alpha-other-bound}
\pucciSub(D^2\phi^a) - \frac{1}{2}y\cdot D\phi^a \leq \left( \frac{(n-1)\lambda}{2\Lambda} + \frac{1}{2} \right) \phi^a \quad \mbox{in} \ \R^n.
\end{equation}
Similar calculations demonstrate that 
\begin{equation}\label{eq:Pucci-bound-alpha-1}
\PucciSub(D^2\phi^a) - \frac{1}{2}y\cdot D\phi^a = a\left(2\Lambda n - 4a\Lambda |y|^2 + |y|^2\right) \phi^a\quad \mbox{for} \  |y| \leq (2a)^{-\frac{1}{2}},
\end{equation}
and 
\begin{equation}\label{eq:Pucci-bound-alpha-2}
\PucciSub(D^2\phi^a) - \frac{1}{2}y\cdot D\phi^a = a\left(2\Lambda (n-1) +2\lambda - 4a\lambda |y|^2 + |y|^2\right) \phi^a \quad \mbox{for} \ |y| \geq (2a)^{-\frac{1}{2}}.
\end{equation}
Thus for any $a \geq (4\lambda)^{-1}$, we have
\begin{equation}\label{eq:Pucci-bound-alpha}
\PucciSub(D^2\phi^a) - \frac{1}{2}y\cdot D\phi^a \leq (2a\Lambda n) \phi^a \quad \mbox{in} \ \R^n.
\end{equation}
Finally, for $a = (4\lambda)^{-1}$, we calculate
\begin{equation}\label{eq:Pucci-alpha-other-bound}
\PucciSub(D^2\phi^a) - \frac{1}{2} y \cdot D\phi^a \geq \left( \frac{(n-1)\Lambda}{2\lambda} + \frac{1}{2} \right) \phi^a \quad \mbox{in} \ \R^n.  \qed
\end{equation}

\medskip

For the rest of this section, we let $\scon$ denote the special constant
\begin{equation*}
\scon := \frac{1}{8\Lambda}.
\end{equation*}
From \EQ{pucci-bound-alpha-1} and \EQ{pucci-bound-alpha-2}, we see that $\phi^b$ satisfies
\begin{equation}\label{eq:pucci-bound-alpha-b}
\pucci(D^2\phi^b) - \frac{1}{2} y\cdot D\phi^b \geq \left( \frac{n\lambda}{4\Lambda} + \frac{1}{16\Lambda} |y|^2 \right) \phi^b \quad \mbox{in} \ \R^n.
\end{equation}
Define a norm $\normx{\cdot}$ on $C(\R^n)$ by
\begin{equation*}
\normx{u} = \sup_{y\in \R^n} \frac{|u(y)|\exp\left(b |y|^2\right)}{1+|y|^2},\end{equation*}
and let $X$ denote the Banach space
\begin{equation*}
X = \left\{ u \in C(\R^n) : \normx{u} < \infty \right\}.
\end{equation*}
Notice that convergence in $X$ implies uniform convergence in $\R^n$. Define the set
\begin{equation*}
\decayseta := \left\{ u \in X : 0 \leq u(y) \leq C \exp( -b |y|^2) \ \mbox{for some} \ C>0 \right\}.
\end{equation*}
Notice that $\decayseta$ is a convex subset of $X$.

\begin{prop}\label{prop:unique-solution-allspace}
For each $v \in X$ such that $v\geq 0$, there exists a unique solution $u \in X$ of the equation
\begin{equation}\label{eq:unique-solution-allspace}
F(D^2u)  -\frac{1}{2} y\cdot Du = v \quad \mbox{in} \ \R^n.
\end{equation}
Moreover, $u \in \decayseta$.
\end{prop}
\proof
We will first demonstrate existence. For each $R> 0$, let $u^R$ be the unique solution of the Dirichlet problem
\begin{equation*}
\left\{ \begin{aligned}
& F\left(D^2u_R\right) - \frac{1}{2}y \cdot Du_R = v & \mbox{in} & \ B_R,\\
& u_R = 0 & \mbox{on} & \ \partial B_R.
\end{aligned} \right.
\end{equation*}
Set  $K:= \normx{v} \max \left\{ 4\Lambda / n\lambda ,16\Lambda \right\}$. According to \EQ{pucci-bound-alpha-b}, the function $\psi:= K\phi^b$ is a supersolution of
\begin{equation*}
F(D^2\psi) - \frac{1}{2} y \cdot D\psi \geq \normx{v} \left( 1 + |y|^2 \right) \phi^b  \geq v \quad \mbox{in} \ \R^n.
\end{equation*}
By the maximum principle, $0 \leq u_R \leq \psi = K \exp(-b|y|^2)$ for every $R>0$. Using local $C^\alpha$ estimates for uniformly elliptic equations (c.f. \cite{Winter:2009}), we deduce that for each fixed $R_0> 0$,
\begin{equation*}
\sup_{ R > R_0 +1} \| u_R \|_{C^\alpha(B_{R_0})} < \infty.
\end{equation*}
Extend $u_R$ to be zero outside $\R^n \backslash B_R$, and extract a subsequence $R_j \to \infty$ such that
\begin{equation*}
u_{R_j} \rightarrow u \quad \mbox{locally uniformly in} \ \R^n
\end{equation*}
for some function $u \in C(\R^n)$. Evidently $0 \leq u \leq \psi$, and thus $u \in \decayseta$. From the stability properties of viscosity solutions under uniform convergence, it follows that $u$ is a solution of equation \EQ{unique-solution-allspace}. Uniqueness follows from \COR{parabolic-HCP}. Indeed, if $u_1, u_2\in \decayseta$ are solutions of \EQ{unique-solution-allspace}, then the function $w:= u_1-u_2$ satisfies
\begin{equation*}
\pucci(D^2w) - \frac{1}{2} y\cdot Dw \leq 0 \quad \mbox{in} \ \R^n.
\end{equation*}
Comparing $w$ with $\phi^{2b}$ and using \COR{parabolic-HCP}, we deduce that $w\leq 0$ in $\R^n$.
\qed

\medskip

We denote by $X_+$ the set
\begin{equation*}
X_+ = \left\{ u \in X : u\geq 0 \right\}.
\end{equation*}
Let $\mathcal{A}:X_+ \to X_+$ be the solution operator of \EQ{unique-solution-allspace}. That is, $\mathcal{A}(v) := u$, where $u$ is the unique solution of \EQ{unique-solution-allspace}. Then $\mathcal{A}\left( X_+ \right) \subseteq \decayseta$. It will be convenient to use the notation
\begin{equation*}
\Falpha{u} := F(D^2u) - \frac{1}{2} y\cdot Du - \alpha u.
\end{equation*}
Define the constant
\begin{equation}\label{eq:def-of-alpha-plus}
\alpha^+(F) := \sup \left\{ \alpha : \mbox{there exists} \ \varphi \in X_+\backslash \{ 0\}  \ \mbox{such that} \ \Falpha{\varphi} \geq 0 \ \mbox{in} \ \R^n \right\}.
\end{equation}
We call $\alpha^+(F)$ the \emph{positive anomalous exponent} of $F$. When there is no ambiguity, we will drop the dependence on $F$ and write $\alpha^+ = \alpha^+(F)$. From \COR{parabolic-HCP} and Hopf's Lemma, it is clear that the anomalous exponent satisfies
\begin{equation*}
\alpha^+ \leq \inf\left\{ \alpha : \mbox{there exists} \ \varphi \in X_+\backslash \{ 0\} \ \mbox{such that} \ \Falpha{\varphi} \leq 0 \ \mbox{in} \ \R^n \right\}.
\end{equation*}
From \LEM{rugged-calculations} we see that
\begin{equation}\label{eq:pucci-alpha-bounds}
\frac{n\lambda}{2\Lambda} \leq \alpha^+\left( \pucciSub \right) \leq \frac{(n-1) \lambda + \Lambda}{2\Lambda} \leq \frac{n}{2} \leq \frac{(n-1) \Lambda + \lambda}{2\lambda} \leq \alpha^+(\PucciSub) \leq \frac{n\Lambda}{2\lambda}.
\end{equation}
Moreover, if $\lambda \neq \Lambda$, then all of the inequalities in \EQ{pucci-alpha-bounds} are strict. If $F$ and $G$ are two uniformly elliptic, positively homogeneous operators such that $F \leq G$, then
\begin{equation*}
\alpha^+(F) \leq \alpha^+(G).
\end{equation*}
Thus from \EQ{pucci-alpha-bounds} we deduce that for every operator $F$ satisfying \EQ{Felliptic} and \EQ{Fhomogeneous},
\begin{equation}
\frac{n\lambda}{2\Lambda} \leq \alpha^+(F) \leq \frac{n\Lambda}{2\lambda}.
\end{equation}

We will show that $\alpha^+$ is an eigenvalue of the operator $\mathcal{A}$, using the Leray-Schauder alternative. For the convenience of the reader, we first state this result. A proof can be found in \cite{Granas:Book}.

\begin{definition}
If $Y$ and $Z$ are Banach spaces, we say a (possibly nonlinear) map $\mathcal{A}:Y \to Z$ is \emph{compact} if, for each bounded subset $B \subseteq Y$, the closure of the set $\{ \mathcal{A}(x) : x \in B \}$ is compact in $Z$.
\end{definition}

\begin{thm}[Leray-Schauder Alternative]\label{thm:leray-schauder-alt}
Suppose $Y$ is a Banach space, and $C \subseteq Y$ is a convex subset of $Y$ such that $0 \in C$. Assume that $\mathcal{A}:C \to C$ is a (possibly nonlinear) function which is compact and continuous. Then at least one of the following holds:
\begin{enumerate}
\item the set $\{ x \in C \, : \, x = \mu \mathcal{A}(x) \mbox{ for some } 0 < \mu < 1 \}$ is unbounded in $Y$,
\end{enumerate}
or
\begin{enumerate}
\addtocounter{enumi}{1}
\item there exists $x\in C$ for which $x = \mathcal{A}(x)$.
\end{enumerate}
\end{thm}

In order to apply \THM{leray-schauder-alt} in our setting, we must verify that the nonlinear operator $\mathcal{A}$ is continuous and compact.

\begin{prop}
The operator $\mathcal{A}$ is continuous and compact with respect to $\normx{\cdot}$.
\end{prop}
\proof
Let $\{ v_k \}_{k\geq 1} \subseteq X$ such that $\normx{v_k} \leq 1$. Set $u_k := \mathcal{A}(v_k)$. Let $\ep > 0$ be given, and fix a large constant $R> 0$ to be selected below. As in the proof of \PROP{unique-solution-allspace}, from \EQ{pucci-bound-alpha-1} and \EQ{pucci-bound-alpha-2} we see that the function $\phi:=M\phi^{b}$ satisfies
\begin{equation*}
\pucci(D^2\phi) - \frac{1}{2}y\cdot D\phi \geq \left( 1 + |y|^2 \right) \phi^{b} \geq v_k \quad \mbox{in} \ \R^n
\end{equation*}
for $M:= \max\left\{ 4\Lambda / n\lambda, 16\Lambda \right\}$ and every $k\geq 1$. It follows that $u_k \leq \psi$ for every $k\geq 1$. Using local $C^\alpha$ estimates, we have
\begin{equation*}
\sup_{k \geq 1} \| u_k \|_{C^\alpha(B_R)} < \infty.
\end{equation*}
Therefore, we may select a subsequence, which we also denote by $k$, such that 
\begin{equation*}
\lim_{K\to \infty} \sup_{k,l \geq K} \| u_k - u_l \|_{L^\infty(B_R)} = 0.
\end{equation*}
Now take $R = \left( 2M/ \ep\right)^{1/2}$. Then for any $y \geq R$ and $k,l \geq 1$,
\begin{equation*}
\frac{\left| u_k(y) - u_l (y)\right| \exp(b|y|^2)}{1+|y|^2} \leq \frac{2|\psi(y)| \exp(b|y|^2)}{1+R^2} = \frac{2M}{1+R^2} \leq \ep.
\end{equation*}
It follows that 
\begin{equation*}
\lim_{K\to \infty} \sup_{k,l \geq K} \normx{ u_k - u_l} \leq \ep.
\end{equation*}
A diagonalizing procedure now produces a subsequence of $\{ u_k \}$ which is Cauchy in $X$. Therefore, $\mathcal{A}$ is compact.

To see that $\mathcal{A}$ is continuous, suppose in addition that the sequence $v_k$ converges strongly in $X$ to a function $v\in X$. In particular, $v_k \rightarrow v$ uniformly in $\R^n$. We can find $u \in X$ and a subsequence $u_{k_j}$ such that $u_{k_j} \rightarrow u$ in $X$, and hence uniformly in $\R^n$. By the stability properties of viscosity solutions with respect to uniform convergence, it follows that $u=\mathcal{A}(v)$. By uniqueness, the full sequence $u_k$ converges to $u$.
\qed

\begin{prop}\label{prop:principal-eigenfunction-allspace}
There exists a unique $\varphi^+ \in X$ such that $\varphi^+(0) = 1$ and 
\begin{equation}\label{eq:eigenvalue-allspace}
F(D^2\varphi^+) - \frac{1}{2} y\cdot D\varphi^+ = \alpha^+ \varphi^+ \quad \mbox{in} \ \R^n.
\end{equation}
Moreover, $\varphi^+ \in \decayseta \cap C^{1,\alpha}_{\mathrm{loc}}(\R^n)$, and $\varphi^+ \in C^{2,\alpha}_{\mathrm{loc}}(\R^n)$ if $F$ is concave or convex.
\end{prop}
\proof
Select $w \in \decayseta$ such that $\normx{w} = 1$. We claim that for any $\ep > 0$,
\begin{equation}\label{eq:existence-claim}
\mbox{if} \  u \in \decayseta\ \mbox{and} \ \alpha\geq 0 \ \mbox{satisfy} \ u = \alpha \mathcal{A}(u+\ep w), \ \mbox{then} \ \alpha \leq \alpha^+.
\end{equation}
Indeed, for such $\alpha \geq 0$ and $u \in \decayseta$, we have
\begin{equation*}
\Falpha{u} \geq  0 \quad \mbox{in} \ \R^n.
\end{equation*}
Since $w\not\equiv 0$, if $\alpha > 0$ then $u \not\equiv 0$. We see from definition \EQ{def-of-alpha-plus} that in this case $\alpha \leq \alpha^+$. Obviously if $\alpha = 0$, then $u\equiv 0$. Our claim \EQ{existence-claim} is confirmed.

\smallskip

We may now apply \THM{leray-schauder-alt} to deduce that for each $\ep > 0$, the set
\begin{equation*}
D_\ep := \left\{ u \in \decayseta : \mbox{there exists} \ 0 \leq \alpha \leq \alpha^+ + \ep \ \mbox{such that} \ u = \alpha \mathcal{A}(u+\ep w) \right\}
\end{equation*}
is unbounded in $X$. Select $u_\ep \in D_\ep$ such that $\normx{u_\ep} \geq 1$. Let $\alpha_\ep\geq 0$ such that 
\begin{equation*}
u_\ep = \alpha_\ep \mathcal{A}(u_\ep+\ep w)
\end{equation*}
Evidently, $\alpha_\ep > 0$. Normalize by setting $v_\ep := u_\ep / \normx{u_\ep}$, and notice that by the homogeneity of $\mathcal{A}$, the function $v_\ep$ satisfies
\begin{equation*}
v_\ep = \alpha_\ep \mathcal{A}\left(v_\ep + \ep w / \normx{u_\ep} \right).
\end{equation*}
By the compactness of $\mathcal{A}$, we may select $\varphi^+ \in X$, a number $0 \leq \alpha^* \leq \alpha^+$, and a subsequence $\ep_j\to 0$, such that 
\begin{equation*}
v_{\ep_j} \to \varphi^+ \ \mbox{in} \ X \quad \mbox{and} \quad \alpha_{\ep_j} \to \alpha^*.
\end{equation*}
Since $\mathcal{A}$ is continuous, it follows that $\varphi^+ = \alpha^* \mathcal{A}(\varphi^+)$. Thus $\varphi^+\in \decayseta$. Clearly $\normx{\varphi^+} = 1$, and thus $\alpha^* > 0$. By Hopf's Lemma, $\varphi^+> 0$ in $\R^n$. 

We will now argue that $\alpha^* = \alpha^+$, and that $\varphi^+$ is unique up to multiplication by a positive constant. Suppose that $\alpha \geq \alpha^*$ and $\psi \in X_+ \backslash \{ 0 \}$ are such that
\begin{equation*}
F(D^2\psi) - \frac{1}{2} y\cdot D\psi \geq \alpha \psi \quad \mbox{in} \ \R^n.
\end{equation*}
Then
\begin{equation}
F_{\alpha^*} \left[ \varphi^+ \right] = 0 \leq \Falpha{\psi} \leq F_{\alpha^*} \left[ \psi \right].
\end{equation}
According to \COR{parabolic-HCP}, we have $\psi \equiv t \varphi^+$ for some $t > 1$. This implies that $\alpha^* = \alpha$. Recalling \EQ{def-of-alpha-plus}, we see that $\alpha^* \geq \alpha^+$. Recalling that by construction $\alpha^* \leq \alpha^+$, we deduce that $\alpha^* = \alpha^+$. Moreover, we have shown that $\varphi^+$ is unique up to multiplication by a positive constant. 

The last statement in the proposition follows from the standard regularity theory for uniformly elliptic equations (c.f. \cite{Caffarelli:Book,Trudinger:1988}).
\qed

\proof[Proof of \THM{existence-self-similar-solutions}]
Define
\begin{equation}
\Phi^+(x,t) := t^{-\alpha^+} \varphi^+ \left( \frac{x}{\sqrt{t}}\right).
\end{equation}
Assuming that $\Phi^+$ and $\varphi^+$ are smooth, and using \EQ{eigenvalue-allspace}, we easily verify that $\Phi^+$ is a solution of \EQ{fully-nonlinear-parabolic}. If $\Phi^+$ and $\varphi^+$ are not smooth, our calculation can be made rigorous in the viscosity sense by the use of smooth test functions. The uniqueness of $\alpha^+$ and $\Phi^+$ is established by performing this computation in reverse and appealing to \PROP{principal-eigenfunction-allspace}. All of the corresponding assertions regarding $\alpha^-$ and $\Phi^-$ now follow from \REM{its-all-the-same}.
\qed

\medskip

We conclude this section with an estimate of our self-similar solution $\Phi^+$ from above and below, and an example.

\begin{lem}\label{lem:bound-above-and-below}
For each $0 < a < (4\Lambda)^{-1}$, there exists a constant $C>0$ such that
\begin{equation}\label{eq:bounded-by-gaussian-above}
\varphi^+(y) \leq C \exp\left( -a|y|^2 \right).
\end{equation}
Likewise, for each $a > (4\lambda)^{-1}$, there exists a constant $C>0$ such that
\begin{equation}\label{eq:bounded-by-gaussian-below}
\exp\left( -a|y|^2\right) \leq C \varphi^+(y).
\end{equation}
\end{lem}
\proof
By construction, since $\varphi^+ \in \decayseta$ we have that the estimate \EQ{bounded-by-gaussian-above} holds for $a_1 = (8\Lambda)^{-1}$. We will therefore only show \EQ{bounded-by-gaussian-below}, as a similar argument obtains \EQ{bounded-by-gaussian-above} for all $a_1 < (4\Lambda)^{-1}$. For $a > (4\lambda)^{-1}$, by \EQ{Pucci-bound-alpha-2} we have that 
\begin{equation*}
\Pucci(D^2 \phi^a) - \frac{1}{2} y\cdot D\phi^a \leq \alpha^+ \phi^a \quad \mbox{in} \ \R^n \backslash B_r,
\end{equation*}
provided that we take $r> 0$ so large that
\begin{equation*}
r^2 > (2a)^{-1} \quad \mbox{and} \quad r^2 \geq \frac{a\left( 2\Lambda(n-1) + 2\lambda\right) - \alpha^+}{a(4a\lambda -1)}.
\end{equation*}
Also take $r>R$, where $R = R\left(\alpha^+,\Lambda\right)$ is the constant in \LEM{SDMP}.
Let $C$ be so large that $\phi^a \leq C \varphi^+$ on $B_r$. Then the function $w:= \phi^a - C\varphi^+$ satisfies
\begin{equation*}
\pucci(D^2w) - \frac{1}{2} y\cdot Dw \leq \alpha^+ w \quad \mbox{in} \ \R^n \backslash B_r,
\end{equation*}
and $w\leq 0$ on $\partial B_r$. According to \LEM{SDMP}, the function $w\leq 0$ in $\R^n$. That is, $\phi^a \leq C \varphi^+$ in $\R^n$. \qed

\begin{cor}\label{cor:bound-above-and-below}
For each $0 < a_1 < (4\Lambda)^{-1} \leq (4\lambda)^{-1} < a_2$, there exists a constant $C>1$ such that
\begin{equation}\label{eq:Phi-bound-above-and-below}
C^{-1} t^{-\alpha^+} \exp \left( -a_2 |x|^2 / t\right) \leq \Phi^+(x,t) \leq C t^{-\alpha^+} \exp \left( -a_1 |x|^2 / t\right) 
\end{equation}
for all $(x,t) \in \R^n\times\R_+$.
\end{cor}

\begin{example}
Consider the case that $F$ is convex. Then $F$ is a supremum of a collection of linear operators $L^k$ with constant coefficients, and each of which satisfy \EQ{Felliptic} and \EQ{Fhomogeneous}. Since $\alpha^+(L^k) = \alpha^-(L^k) = n/2$ for every $k$, we deduce that
\begin{equation*}
\alpha^-(F) \leq \frac{n}{2} \leq \alpha^+(F).
\end{equation*}
We claim that these inequalities are strict unless $F$ is linear. Suppose that $\alpha^+(F) = n/2$. 
Let $\varphi$ and $\varphi_k$ be the functions obtained in \PROP{principal-eigenfunction-allspace} for $F$ and $L^k$, respectively. Notice that
\begin{equation*}
F_{n/2} \left[\varphi\right] = 0 = L^k_{n/2} \left[\varphi_k\right] \leq F_{n/2} \left[\varphi_k\right] \quad \mbox{in} \ \R^n. 
\end{equation*}
According to \COR{parabolic-HCP}, $\varphi \equiv \varphi_k$ for every $k$. That is, the fundamental solutions of the constant-coefficient linear parabolic operators $L^k$ are equal. This implies that $L^k=L$ for every $k$ (see Friedman \cite{Friedman:Book:1964}). Hence $F=L$.
\end{example}

\begin{remark}
In the case that $F(M)$ depends only on the eigenvalues of $M$, we immediately deduce that $\varphi^+$ and hence $\Phi^+(\cdot, t)$ are radial functions. This follows from the invariance of the equation under an orthogonal change of variables, and \COR{parabolic-HCP}. In particular, the self-similar solutions corresponding to the operators $\Pucci$ and $\pucci$ are radial.
\end{remark}


\section{Asymptotic convergence to self-similar solutions} \label{sec:asymptotics}

In this section, we present the proof of \THM{convergence-self-similar-solutions}. Owing to \REM{its-all-the-same}, we need only prove the first statement. For ease of notation, we write $\alpha = \alpha^+(F)$ and $\Phi = \Phi^+$. Fix a solution $u=u(x,t)$ of the equation
\begin{equation}\label{eq:PDE}
u_t + F(D^2u) = 0 \quad \mbox{in} \quad \R^n\times\R_+, \\
\end{equation}
subject to the initial condition
\begin{equation}
u(x,0) = g(x).
\end{equation}
We require the initial data $g$ to be continuous, not identically zero, and to satisfy the condition
\begin{equation*}
0 \leq g(x) \leq C_0 e^{-B|x|^2}
\end{equation*}
for some constants $B,C_0>0$. For $\sigma > 0$, we denote 
\begin{equation*}
u^\sigma: = \rescale{\sigma} u (x,t) := \sigma^{\alpha} u\left( \sigma^{1/2} x , \sigma t\right).
\end{equation*}
For each $\sigma > 0$, the function $u^\sigma$ is a solution of \EQ{PDE}.

\medskip

We intend to show that as the parameter $\sigma \to \infty$, the rescaled solutions $u^\sigma$ converge locally uniformly in $\R^n\times\R_+$ to a positive multiple of $\Phi(x,t)$. Recall that $\Phi$ is invariant under $\rescale{\sigma}$:
\begin{equation}\label{eq:Phi-invariant}
\Phi(x,t) = \rescale{\sigma} \Phi(x,t) = \sigma^{\alpha} \Phi\left(\sigma^{1/2} x, \sigma t \right) \quad \mbox{for all} \quad (x,t) \in \R^n\times\R_+, \  \sigma > 0.
\end{equation}

\medskip

The proof of \THM{convergence-self-similar-solutions} will consist of a series of lemmas. As a preliminary step, we show that $u$ is bounded between positive multiples of $\Phi$, possibly shifted in time.

\begin{lem}\label{lem:bound-above-initial}
For each $\tau > \frac{1}{4\lambda B}$ there exists a constant $C > 0$, depending only on $C_0$ and $\tau$, such that 
\begin{equation}\label{eq:bound-above-initial}
u(x,t) \leq C \Phi(x,t + \tau) \quad \mbox{for all} \quad (x,t) \in \R^n\times\R_+.
\end{equation}
\end{lem}
\proof
If $\tau > \frac{1}{4B\lambda}$, then according to \COR{bound-above-and-below},
\begin{equation*}
\Phi(x,\tau) \geq c(\tau) e^{-B|x|^2},
\end{equation*}
provided we choose $c(\tau) > 0$ small enough. Thus $C\Phi(x,\tau) \geq g(x)$ for $C:= C_0 / c(\tau)$ and $x\in \R^n$. The maximum principle implies that $C\Phi (x,t+\tau) \geq u(x,t)$ for all $(x,t) \in \R^n\times\R_+$.
\qed

\begin{lem}\label{lem:bound-below-initial}
For each $t_0,\tau > 0$, there exists $C > 0$, depending only on $B$, $C_0$, $t_0$, and $\tau$, such that
\begin{equation}\label{eq:bound-below-initial}
\Phi(x,t) \leq C u(x,t+\tau)  \quad \mbox{for all} \quad x\in \R^n, \ t\geq t_0.
\end{equation}
\end{lem}
\proof
By the strong maximum principle, $u(x,\tau) > 0$ on $\R^n$. Let $C>0$ be so large that
\begin{gather*}
\Phi(x,t_0) \leq C u(x,t_0+\tau) \quad \mbox{for all} \ |x| \leq 1, \intertext{and}
\Phi(x,t) \leq C u(x,t+\tau) \quad \mbox{for all} \ |x| = 1, \ 0 < t \leq t_0.
\end{gather*}
Applying the maximum principle, we have $\Phi(x,t_0) \leq C u(x,t_0 + \tau)$ for all $x \in \R^n$, and \EQ{bound-below-initial} follows from another application of the maximum principle.
\qed

\medskip

For $\tau = 1/(2\lambda B)$, we use \EQ{Phi-invariant} to rewrite the inequality \EQ{bound-above-initial} in terms of $u^\sigma$ as
\begin{equation}\label{eq:u-sigma-bounded-above-by-Phi}
u^\sigma (x,t) \leq C \Phi(x, t+ \tau/ \sigma) \quad \mbox{for all} \quad (x,t) \in \R^n\times\R_+, \ \sigma > 0.
\end{equation} 
Recalling \EQ{Phi-bound-above-and-below}, we see that for some constant $C(t) > 0$ depending only on a lower bound for $t>0$, in addition to $B$ and $C_0$, we have the estimate
\begin{equation}\label{eq:u-sigma-bounded-above}
u^\sigma(x,t) \leq C(t) \exp \left( - |x|^2 / 8\Lambda ( t + \tau / \sigma) \right) \quad \mbox{for all} \quad (x,t) \in \R^n\times\R_+, \ \sigma \geq 1.
\end{equation}
According to \EQ{u-sigma-bounded-above} and local H\"older estimates for solutions of uniformly parabolic equations (see Wang \cite[Theorem 4.19]{Wang:1992a}), we obtain
\begin{equation*}
\sup_{\sigma \geq 1} \| u^\sigma \|_{C^\gamma( \bar{Q} ) } < \infty
\end{equation*}
for some $0 < \gamma < 1$ and any compact parabolic domain $\bar{Q} \subseteq \R^n\times\R_+$. Therefore, for every sequence $\sigma_k \to \infty$, we may select a function $U\in C(\R^n\times\R_+)$ and a subsequence, also denoted by $\sigma_k$, such that $u^{\sigma_k} \to U$ locally uniformly in $\R^n\times\R_+$. By the stability of viscosity solutions with respect to uniform convergence, each such rescaled limit $U$ is a solution of equation \EQ{PDE}. 

\medskip

Let $\limitset$ denote the set of such sequential limits $\{ U \}$ of the family $\{ u^\sigma \}_{\sigma \geq 1}$. We will prove \THM{convergence-self-similar-solutions} by showing that $\mathcal{S}$ is a singleton set consisting only of a positive multiple of $\Phi$.

\begin{lem}\label{lem:convergence-squeeze}
There exists a positive constants $C>0$ such that for all $U \in \mathcal{S}$,
\begin{equation}\label{eq:convergence-squeeze}
C^{-1}\Phi \leq U \leq C \Phi \quad \mbox{in} \ \R^n\times\R_+.
\end{equation}
\end{lem}
\proof
From \EQ{u-sigma-bounded-above-by-Phi}, we see that for each $U \in \limitset$ and $(x,t) \in \R^n\times\R_+$,
\begin{equation*}
U(x,t) \leq \limsup_{\sigma \to \infty} u^\sigma ( x, t ) \leq C \Phi\left( x,t\right).
\end{equation*}
For the other direction, fix $t_0, \tau > 0$. According to \EQ{bound-below-initial}, for each $t > 0$ the inequality 
\begin{equation}\label{eq:convergence-squeeze-1}
\sigma^{-\alpha} \Phi\left( \sigma^{1/2} x, \sigma t \right) \leq C \sigma^{-\alpha} u\left( \sigma^{1/2} x , \sigma t + \tau \right)
\end{equation}
holds for all $(x,t) \in \R^n\times\R_+$ and all sufficiently large $\sigma \geq 1$. Rewrite \EQ{convergence-squeeze-1} as 
\begin{equation*}
\Phi(x,t) \leq C u^\sigma \left( x, t + \tau / \sigma\right).
\end{equation*}
Take the lim-inf of the right side as $\sigma \to \infty$ to see that $\Phi(x,t) \leq C U( x, t)$ for any $U \in \limitset$.
\qed

\medskip

Notice that \LEM{convergence-squeeze} and local H\"older estimates imply that 
\begin{equation}\label{eq:limitset-estimate}
\sup_{U \in \limitset} \left\| U \right\|_{C^\gamma(\bar{Q})} < \infty,
\end{equation}
for every compact subset $\bar{Q} \subseteq \R^n\times\R_+$.

\medskip

Define the constant
\begin{equation}\label{eq:Cstar}
C^* : = \inf \left\{ C>0 : \mbox{there exists} \ U \in \limitset \ \mbox{such that} \ U \leq C \Phi \right\}.
\end{equation}
In light of \LEM{convergence-squeeze}, $0 < C^* < \infty$. We will eventually show that $\limitset = \left\{ C^* \Phi\right\}$. There are two basic steps in the proof. First, we will show that $U \leq C^* \Phi$ for every $U \in \limitset$. Second, we show that if $U \not\equiv C^* \Phi$ for some $U \in \limitset$, then we can find another function $V \in \limitset$ and a small number $\delta > 0$ such that $V \leq (C^*- \delta) \Phi$, in contradiction to the definition \EQ{Cstar} of $C^*$. Most of the subtlety in the proofs of these statements arise from difficulties in managing the ``tails" of $\Phi$. These obstructions are removed by the construction of a special subsolution, which we use as a comparison function.

\begin{lem}\label{lem:special-subsolution}
For each $a > \frac{1}{4\lambda}$, there exists $r,\eta >0$ and a subsolution $w = w(y,s)$ of the differential inequality
\begin{equation} \label{eq:special-subsolution}
w_s + \PucciSub(D^2w) - \frac{1}{2} y\cdot Dw \leq 0 \quad \mbox{in} \ \R^n\times\R_+,
\end{equation}
satisfying the initial conditions
\begin{equation*}
w(y,0) \leq e^{-a|y|^2} \quad \mbox{for all} \ |y| \leq r, \quad \mbox{and} \quad w(y,0) \leq -\eta e^{-|y|} \quad \mbox{for all} \ |y| > r,
\end{equation*}
and such that for each $R>0$ there exists a time $S>0$ such that 
\begin{equation*}
w(y,s) > 0 \quad \mbox{for every} \ |y| \leq R, \ s\geq S.
\end{equation*}
\end{lem}
\proof
Select a constant $a > (4\lambda)^{-1}$ and set $\varphi(y) := \phi^a (y) = e^{-a|y|^2}$. Recall from \EQ{Pucci-bound-alpha} that
\begin{equation*}
\Pucci ( D^2\varphi) - \frac{1}{2} y \cdot D\varphi \leq \left( 2a\Lambda n \right) \varphi \quad \mbox{in} \ \R^n.
\end{equation*}
Let $\beta:= 1+ 2a\Lambda n$ and $r_1:= 2(\beta+ \Lambda+1)$, and define a function $\psi(y) = \min\left\{ e^{-r_1} , e^{-|y|} \right\}$. Recalling \EQ{SDMP-1} and that the minimum of supersolutions is a supersolution in the viscosity sense, we see that
\begin{equation*}
\pucci(D^2\psi) - \frac{1}{2}y \cdot D\psi \geq 0 \quad \mbox{in} \ B_{r_1+1},
\end{equation*}
and
\begin{equation*}
\pucci(D^2\psi ) - \frac{1}{2} y \cdot D\psi \geq (\beta + 1) \psi \quad \mbox{in} \ \R^n \backslash B_{r_1}.
\end{equation*}
Now define $\bar{\varphi}(y,s) := e^{-\beta s} \varphi(y)$ and $\bar{\psi}(y,s) := - e^{-(\beta+1)s} \psi(y)$. Then
\begin{equation*}
\bar{\varphi}_s + \Pucci(D^2\bar{\varphi}) - \frac{1}{2} y\cdot D\bar{\varphi} \leq - \bar{\varphi} \quad \mbox{in} \ \R^n\times\R_+,
\end{equation*}
and 
\begin{equation*}
\bar{\psi}_s + \Pucci(D^2\bar{\psi}) - \frac{1}{2} y\cdot D\bar{\psi} \leq (\beta+1)e^{-(\beta+1)s}e^{-r_1} \chi \quad \mbox{in} \ \R^n\times\R_+,
\end{equation*}
where $\chi\equiv 1$ on $\bar{B}_{r_1}$ and $\chi\equiv 0$ on $\R^n \backslash \bar{B}_{r_1}$. Set
\begin{equation*}
\delta := \frac{1}{\beta+1} e^{r_1 - ar_1^2} > 0.
\end{equation*}
Then $\bar{\varphi}(y,s) \geq \delta (\beta+1) e^{-(\beta+1)s -r_1}$ for all $y\in \bar{B}_{r_1}$ and $s\geq 0$. Therefore, the function $w:= \bar{\varphi} + \delta \bar{\psi}$ satisfies \EQ{special-subsolution}.

We now investigate the set of $(y,s)$ for which $w > 0$. For every $|y| > r_1$, 
\begin{equation*}
w(y,s) = e^{-\beta s} \left( e^{-a|y|^2} - \delta e^{-s} e^{-|y|} \right).
\end{equation*}
From this expression, we observe that $w(y,s) > 0$ whenever $s > a|y|^2 - |y| + \log \delta$ and $|y|>r_1$. Finally, select $r >0$ large enough that $r \geq r_1$ and $as^2 \geq s + \log\frac{2}{\delta}$ whenever $s \geq r$. This choice of $r$ ensures that
\begin{equation*}
w(y,0) \leq -\frac{\delta}{2} e^{-|y|} \quad \mbox{for all} \  |y| > r.
\end{equation*}
Taking $\eta := \delta /2$, the proof is complete.
\qed

\begin{cor}\label{cor:magic}
There exist $r,\eta > 0$ such that for any $R>0$ and any subsolution $u$ of
\begin{equation*}
u_t + \pucci(D^2u) \leq 0 \quad \mbox{in} \ \R^n \times (1,\infty)
\end{equation*}
satisfying initial conditions
\begin{equation*}
u(x,1) \leq -1 \quad \mbox{for every} \ |x| \leq r, \quad \mbox{and} \quad u(x,1) \leq \eta e^{-|x|} \quad \mbox{for every} \ |x| > r,
\end{equation*}
there exists $T> 1$ such that
\begin{equation*}
u(x,t) < 0 \quad \mbox{for all} \quad t\geq T, \ |x| \leq R \sqrt{t}.
\end{equation*}
\end{cor}
\proof
Let $r$, $\eta$, and $w$ be as in \LEM{special-subsolution} for $a = \frac{1}{2\lambda}$. Define
\begin{equation*}
v(x,t) := - w\left( \frac{x}{\sqrt{t}}, \log t\right), \quad (x,t) \in \R^n \times (1,\infty).
\end{equation*}
Then $v$ is a supersolution of the equation
\begin{equation*}
v_t + \pucci(D^2v) \geq 0 \quad \mbox{in} \ \R^n \times (1,\infty),
\end{equation*}
such that
\begin{equation*}
v(x,1) \geq \eta e^{-|x|} \quad \mbox{for every}\ |x| > r, \quad \mbox{and} \quad v(x,1) \geq - e^{-a|x|^2} \quad \mbox{for every}\ |x| \leq r.
\end{equation*}
In particular, $v \geq u$ at time $t=1$ and thus $v\geq u$ in $\R^n \times (1,\infty)$ by the maximum principle. According to the conclusion of \LEM{special-subsolution}, for each $R > 0$ there exists $S > 0$ such that
\begin{equation*}
v\left( x , t \right) = - w\left( t^{-1/2} x, \log t \right) < 0 \quad \mbox{for all} \quad t^{-1/2} | x | \leq R, \ \log t \geq S.
\end{equation*}
Hence the conclusion is obtained for $T = \exp(S)$.
\qed

\begin{lem}\label{lem:Cstar-upperbound}
For any $U\in \limitset$, 
\begin{equation}
U (x,t) \leq C^* \Phi(x,t) \quad \mbox{for all} \quad (x,t) \in \R^n\times\R_+.
\end{equation}
\end{lem}
\proof
Let $r,\eta$ be as in \COR{magic} and fix a small number $\ep > 0$. Recalling \EQ{u-sigma-bounded-above}, we may choose constants $C_1,a>0$ such that
\begin{equation*}
u^{\sigma} (x,1) \leq C_1 \exp(-a|x|^2) \quad \mbox{for all} \ x\in \R^n, \ \sigma \geq 1.
\end{equation*}
Set $m:= \min_{|x|\leq r} \Phi(x,1)$ and select $r_1\geq r$ such that 
\begin{equation*}
\frac{C_1}{m\ep}  e^{-a|x|^2} \leq  \eta e^{-|x|} \quad \mbox{for all} \quad |x| \geq r_1.
\end{equation*}
According to the definition \EQ{Cstar} of $C^*$, we may select $\sigma_1 \geq 1$ such that
\begin{equation*}
u^{\sigma_1} (x,1) \leq  (C^* + \ep) \Phi(x,1) \quad \mbox{for all} \quad |x| \leq r_1.
\end{equation*}
Define $w:= u^{\sigma_1} - (C^* + 2\ep) \Phi$. Then $w$ is a subsolution of the parabolic equation
\begin{equation*}
w_t + \pucci(D^2w) \leq 0 \quad \mbox{in} \ \R^n \times (1,\infty),
\end{equation*}
and at time $t=1$ the function $w$ satisfies
\begin{gather*}
w(x,1) \leq - \ep \Phi(x,1) \leq  - m \ep \quad \mbox{for all} \quad |x| \leq r, \\
w(x,1) \leq -\ep \Phi(x,1) \leq 0 \quad \mbox{for all} \quad |x| \leq r_1,
\intertext{and}
w(x,1) \leq u^{\sigma_1}(x,1) \leq (m\ep) \eta e^{-|x|} \quad \mbox{for all} \quad |x| > r_1.
\end{gather*}
According to \COR{magic}, for each $R>0$ there exists a time $T=T(R)> 1$ such that 
\begin{equation*}
w(x, t) \leq 0 \quad \mbox{for all} \ |x| \leq Rt^{1/2} \ \mbox{and} \ t\geq T. 
\end{equation*}
This reads
\begin{equation*}
\sigma_1^{\alpha}\, u\left(\sigma_1^{1/2} x, \sigma_1 t\right) \leq (C^* + 2\ep) \Phi(x,t) \quad \mbox{for all} \ |x| \leq Rt^{1/2} \ \mbox{and} \ t\geq T.
\end{equation*}
Thus for any $\sigma > \sigma_1$,
\begin{align*}
u^\sigma(x,t) & = \sigma^\alpha u\left(\sigma^{1/2} x , \sigma t\right) \\
& = \sigma^\alpha u\left( \sigma_1^{1/2} (\sigma/ \sigma_1)^{1/2} x , \sigma_1 (\sigma/\sigma_1) t \right) \\
& \leq (\sigma/\sigma_1)^\alpha \left( C^* + 2\ep \right) \Phi\left( (\sigma/ \sigma_1)^{1/2} x , (\sigma/\sigma_1) t  \right)  \\ 
& = \left( C^* + 2\ep \right) \Phi\left( x,t \right)
\end{align*}
provided that $(\sigma/\sigma_1)^{1/2} |x| \leq R \left( \sigma t /\sigma_1 \right)^{1/2}$ and $\sigma t / \sigma_1 \geq T(R)$. In particular, for each $(x,t) \in \R^n\times\R_+$, there exists $\sigma_2 > \sigma_1$ large enough that
\begin{equation*}
u^\sigma(x,t) \leq \left( C^* + 2\ep \right) \Phi\left( x,t \right) \quad \mbox{for all} \quad \sigma \geq \sigma_2.
\end{equation*}
It follows that for any $U \in \limitset$,
\begin{equation*}
U(x,t) \leq \left( C^* + 2\ep \right) \Phi\left( x,t \right) \quad \mbox{for all} \quad (x,t) \in \R^n\times\R_+.
\end{equation*}
The conclusion is obtained by sending $\ep \to 0$.
\qed

\medskip

To complete the proof of \THM{convergence-self-similar-solutions}, we will need elementary properties of $\limitset$ contained in the following two lemmas.

\begin{lem}\label{lem:limitset-invariant}
If $U \in \limitset$, then $\rescale{\sigma} U \in \limitset$ for any $\sigma > 0$. 
\end{lem}
\proof
Select a sequence $\sigma_j \to \infty$ such that $u^{\sigma_j} \rightarrow U$ locally uniformly in $\R^n\times\R_+$. It is easy to check that for $\tilde{\sigma}_j := \sigma \sigma_j$, the sequence $u^{\tilde{\sigma}_j} \to \rescale{\sigma} U$ locally uniformly in $\R^n\times\R_+$.
\qed

\begin{lem}\label{lem:limitset-closed}
The set $\limitset$ is closed in the topology of local uniform convergence.
\end{lem}
\proof
Let $U_j \in \limitset$ such that $U_j \rightarrow U$ locally uniformly in $\R^n\times\R_+$. Fix a compact subset $\bar{Q}$ of $\R^n\times\R_+$, and select $\sigma_j \to \infty$ such that
\begin{equation*}
\sup_{\bar{Q}} \left| u^{\sigma_j} - U_j \right| \leq 2^{-j}.
\end{equation*}
It is clear that $u^{\sigma_j}$ converges to $U$ as $j\to \infty$, uniformly on $\bar{Q}$. Now a diagonalization argument produces a sequence $\tilde{\sigma}_j \to \infty$ for which the functions $u^{\tilde{\sigma}_j}$ converge to $U$ locally uniformly in $\R^n\times \R_+$, as $j\to \infty$.
\qed

\begin{lem}\label{lem:Cstar-lowerbound}
Suppose that $U \in \limitset$ and $C>0$ are such that $U \leq C \Phi$. Then either $U \equiv C\Phi$ or there exists $\delta > 0$ and $V \in \limitset$ such that $V \leq (C-\delta) \Phi$.
\end{lem}
\proof
Suppose that $U\in \limitset$ and $C> 0$ are such that $U \leq C \Phi$, but $U \not\equiv C\Phi$ in $\R^n\times\R_+$. By the strong maximum principle, $U(x,t) < C\Phi(x,t)$ for every $(x,t) \in \R^n\times\R_+$. Let $r,\eta> 0$ be as in \LEM{special-subsolution}, and choose $\ep > 0$ so small that
\begin{equation*}
U(x,1) \leq (C-\ep) \Phi(x,1) \quad \mbox{for every} \ |x|<r.
\end{equation*}
Let $m:= \min_{|x|\leq r} \Phi(x,1)$. Select $\delta > 0$ small enough that
\begin{equation*}
\frac{\ep - \delta}{\delta^{1/2}} > \frac{\eta}{m} \quad \mbox{\and} \quad \delta^{1/2} \Phi(x,1) \leq e^{-|x|} \ \mbox{for all} \ |x| > r.
\end{equation*}
Denote by $w$ the function
\begin{equation*}
w(x,t) := \frac{1}{\delta^{1/2}} \left( U(x,t) - (C-\delta) \Phi(x,t) \right),
\end{equation*}
which is a subsolution of the equation
\begin{equation*}
w_t + \pucci(D^2w) \leq 0 \quad \mbox{in} \ \R^n\times \R_+.
\end{equation*}
Moreover,
\begin{equation*}
w(x,1) \leq \delta^{1/2} \Phi(x,1) \leq e^{-|x|} \quad \mbox{for every} \ |x| > r,
\end{equation*}
and
\begin{equation*}
w(x,1) \leq -\frac{(\ep-\delta)}{\delta^{1/2}}\Phi(x,t) \leq -\eta \quad \mbox{for every} \ |x| \leq r.
\end{equation*}
According to \COR{magic}, for any $R> 0$ there exists $T(R)> 1$ such that 
\begin{equation*}
U(x,t) \leq \left( C - \delta \right) \Phi(x,t) \quad \mbox{provided that} \ |x|\leq R\sqrt{t}, \ \mbox{and} \ t\geq T.
\end{equation*}
It follows that for each $(x,t) \in \R^n\times\R_+$ there exists $\sigma' > 1$ large enough that
\begin{equation*}
\rescale{\sigma} U (x,t) \leq \left(C - \delta \right) \Phi(x,t) \quad \mbox{for all}\  \sigma \geq \sigma'.
\end{equation*}
According to \LEM{limitset-invariant}, $\rescale{\sigma}(U) \in \limitset$. Recalling \EQ{limitset-estimate}, we may select $V \in C(\R^n\times\R_+)$ such that up to a subsequence, $\rescale{\sigma} U\rightarrow V$ locally uniformly in $\R^n\times\R_+$. It is clear that
\begin{equation*}
V \leq \left( C - \delta \right) \Phi.
\end{equation*}
According to \LEM{limitset-closed}, $V\in \limitset$. 
\qed

\proof[Proof of \THM{convergence-self-similar-solutions}]
According to Lemmas \ref{lem:Cstar-upperbound} and \ref{lem:Cstar-lowerbound} and the definition \EQ{Cstar} of the constant $C^*$, the function $C^* \Phi$ is the only element of $\limitset$. The proof of \THM{convergence-self-similar-solutions} is complete.
\qed

\begin{remark}
Let us repeat a remark made in \cite{Kamin:1991}. If we express the constant $C^*$ obtained in \THM{convergence-self-similar-solutions} as a function of the initial data, $C^* = C^*\!\left[ g \right]$, we see immediately that a nonnegative solution $u$ of \EQ{fully-nonlinear-parabolic} has the property that $t\mapsto C^*\!\left[ u(\cdot , t) \right]$ is constant. We thereby deduce a conservation law for our fully nonlinear equation, generalizing the conservation of mass in the case of a linear operator. It would be interesting to discover more information about $C^*\!\left[ g \right]$ in the general nonlinear case. What is this conserved quantity?

\end{remark}


\section{Acknowledgements}

The authors would like to express their appreciation to their thesis advisor, Lawrence C. Evans for his advice and guidance, and to thank the Department of Mathematics of UC Berkeley, for its support. We also thank Juan Luis V\'azquez for his valuable comments and references, and Grigory Barenblatt for helpful comments. We are also indebted to an anonymous referee whose helpful comments greatly improved this article.
\bibliographystyle{plain}
\bibliography{selfsimilarbib}

\def\cprime{$'$} \def\polhk#1{\setbox0=\hbox{#1}{\ooalign{\hidewidth
  \lower1.5ex\hbox{`}\hidewidth\crcr\unhbox0}}}
\begin{thebibliography}{10}

\bibitem{Armstrong:2009}
S.N. Armstrong.
\newblock Principal eigenvalues and an anti-maximum principle for homogeneous
  fully nonlinear elliptic equations.
\newblock {\em J. Differential Equations}, 246:2958--2987, 2009.

\bibitem{Aronson:1995}
D.~G. Aronson and J.~L. V{\'a}zquez.
\newblock Anomalous exponents in nonlinear diffusion.
\newblock {\em J. Nonlinear Sci.}, 5(1):29--56, 1995.

\bibitem{Barenblatt:Book:1990}
G.~I. Barenblatt, V.~M. Entov, and V.~M. Ryzhik.
\newblock {\em Theory of fluid flows through natural rocks}.
\newblock Kluwer Academic Publishers, Dordrecht, 1990.

\bibitem{Berestycki:1994}
H.~Berestycki, L.~Nirenberg, and S.~R.~S. Varadhan.
\newblock The principal eigenvalue and maximum principle for second-order
  elliptic operators in general domains.
\newblock {\em Comm. Pure Appl. Math.}, 47(1):47--92, 1994.

\bibitem{Birindelli:2006}
I.~Birindelli and F.~Demengel.
\newblock First eigenvalue and maximum principle for fully nonlinear singular
  operators.
\newblock {\em Adv. Differential Equations}, 11(1):91--119, 2006.

\bibitem{Birindelli:2007}
I.~Birindelli and F.~Demengel.
\newblock Eigenvalue, maximum principle and regularity for fully non linear
  homogeneous operators.
\newblock {\em Commun. Pure Appl. Anal.}, 6(2):335--366, 2007.

\bibitem{Caffarelli:Book}
L.~A. Caffarelli and X.~Cabr{\'e}.
\newblock {\em Fully nonlinear elliptic equations}, volume~43 of {\em American
  Mathematical Society Colloquium Publications}.
\newblock American Mathematical Society, Providence, RI, 1995.

\bibitem{Caffarelli:2008}
L.~A. Caffarelli and U.~Stefanelli.
\newblock A counterexample to {$C\sp {2,1}$} regularity for parabolic fully
  nonlinear equations.
\newblock {\em Comm. Partial Differential Equations}, 33(7-9):1216--1234, 2008.

\bibitem{UsersGuide}
M.~G. Crandall, H.~Ishii, and P.-L. Lions.
\newblock User's guide to viscosity solutions of second order partial
  differential equations.
\newblock {\em Bull. Amer. Math. Soc. (N.S.)}, 27(1):1--67, 1992.

\bibitem{Crandall:2000}
M.~G. Crandall, M.~Kocan, and A.~{\'S}wi{\polhk{e}}ch.
\newblock {$L\sp p$}-theory for fully nonlinear uniformly parabolic equations.
\newblock {\em Comm. Partial Differential Equations}, 25(11-12):1997--2053,
  2000.

\bibitem{Friedman:Book:1964}
Avner Friedman.
\newblock {\em Partial differential equations of parabolic type}.
\newblock Prentice-Hall Inc., Englewood Cliffs, N.J., 1964.

\bibitem{Goldenfeld:1990}
N.~Goldenfeld, O.~Martin, Y.~Oono, and F.~Liu.
\newblock Anomalous dimensions and the renormalization group in a nonlinear
  diffusion process.
\newblock {\em Phys. Rev. Lett.}, 64(12):1361--1364, Mar 1990.

\bibitem{Granas:Book}
A.~Granas and J.~Dugundji.
\newblock {\em Fixed point theory}.
\newblock Springer Monographs in Mathematics. Springer-Verlag, New York, 2003.

\bibitem{Haraux:1982}
A.~Haraux and F.~B. Weissler.
\newblock Nonuniqueness for a semilinear initial value problem.
\newblock {\em Indiana Univ. Math. J.}, 31(2):167--189, 1982.

\bibitem{Hulshof:1996}
Josephus Hulshof and Juan~Luis Vazquez.
\newblock Maximal viscosity solutions of the modified porous medium equation
  and their asymptotic behaviour.
\newblock {\em European J. Appl. Math.}, 7(5):453--471, 1996.

\bibitem{Ishii:preprint}
H.~Ishii and Y.~Yoshimura.
\newblock Demi-eigenvalues for uniformly elliptic {I}saacs operators.
\newblock preprint.

\bibitem{Kamin:1991}
S.~Kamin, L.~A. Peletier, and J.~L. V{\'a}zquez.
\newblock On the {B}arenblatt equation of elastoplastic filtration.
\newblock {\em Indiana Univ. Math. J.}, 40(4):1333--1362, 1991.

\bibitem{Kamin:1988}
S.~Kamin and J.~L. V{\'a}zquez.
\newblock Fundamental solutions and asymptotic behaviour for the
  {$p$}-{L}aplacian equation.
\newblock {\em Rev. Mat. Iberoamericana}, 4(2):339--354, 1988.

\bibitem{Lions:1983d}
P.-L. Lions.
\newblock Bifurcation and optimal stochastic control.
\newblock {\em Nonlinear Anal.}, 7(2):177--207, 1983.

\bibitem{Quaas:preprint}
R.~Meneses and A.~Quaas.
\newblock Fujita type exponent for fully nonlinear parabolic equations and
  existence results.
\newblock preprint.

\bibitem{Peletier:1986}
L.~A. Peletier, D.~Terman, and F.~B. Weissler.
\newblock On the equation {$\Delta u+{\frac{1}{2}}x\cdot\nabla u+f(u)=0$}.
\newblock {\em Arch. Rational Mech. Anal.}, 94(1):83--99, 1986.

\bibitem{Quaas:2008}
A.~Quaas and B.~Sirakov.
\newblock Principal eigenvalues and the {D}irichlet problem for fully nonlinear
  elliptic operators.
\newblock {\em Adv. Math.}, 218(1):105--135, 2008.

\bibitem{Trudinger:1988}
N.~S. Trudinger.
\newblock H\"older gradient estimates for fully nonlinear elliptic equations.
\newblock {\em Proc. Roy. Soc. Edinburgh Sect. A}, 108(1-2):57--65, 1988.

\bibitem{Vazquez:Book}
J.~L. V{\'a}zquez.
\newblock {\em The porous medium equation}.
\newblock Oxford Mathematical Monographs. The Clarendon Press Oxford University
  Press, Oxford, 2007.

\bibitem{Wang:1992a}
L.~Wang.
\newblock On the regularity theory of fully nonlinear parabolic equations. {I}.
\newblock {\em Comm. Pure Appl. Math.}, 45(1):27--76, 1992.

\bibitem{Winter:2009}
N.~Winter.
\newblock ${W}^{2,p}$ and ${W}^{1,p}$ estimates at the boundary for solutions
  of fully nonlinear, uniformly elliptic equations.
\newblock {\em Z. Anal. Anwend.}, 28(2):129--164, 2009.

\end{thebibliography}

\end{document}